\documentclass[a4paper,11pt]{article}
\usepackage{amsfonts,amssymb,amsmath}

\usepackage[english]{babel}
 \makeatletter
 \@addtoreset{equation}{section}
 \makeatother

 \makeatletter
 \@addtoreset{enunciato}{section}
 \makeatother

 \newcounter{enunciato}[section]

 \newtheorem{ittheorem}{Theorem}
 \newtheorem{itlemma}{Lemma}
 \newtheorem{itproposition}{Proposition}
 \newtheorem{itcorollary}{Corollary}
 \newtheorem{itdefinition}{Definition}
 \newtheorem{itremark}{Remark}
 \newtheorem{itclaim}{Claim}
 \newtheorem{itfact}{Fact}
 \newtheorem{itconjecture}{Conjecture}

 \newenvironment{theorem}{\addtocounter{enunciato}{1}
 \begin{ittheorem}}{\end{ittheorem}}

 \newenvironment{lemma}{\addtocounter{enunciato}{1}
 \begin{itlemma}}{\end{itlemma}}

 \newenvironment{proposition}{\addtocounter{enunciato}{1}
 \begin{itproposition}}{\end{itproposition}}

 \newenvironment{corollary}{\addtocounter{enunciato}{1}
 \begin{itcorollary}}{\end{itcorollary}}

 \newenvironment{definition}{\addtocounter{enunciato}{1}
 \begin{itdefinition}}{\end{itdefinition}}

 \newenvironment{remark}{\addtocounter{enunciato}{1}
 \begin{itremark}}{\end{itremark}}

 \newenvironment{claim}{\addtocounter{enunciato}{1}
 \begin{itclaim}}{\end{itclaim}}

 \newenvironment{fact}{\addtocounter{enunciato}{1}
 \begin{itfact}}{\end{itfact}}

 \newenvironment{conjecture}{\addtocounter{enunciato}{1}
 \begin{itconjecture}}{\end{itconjecture}}

 \newcommand{\be}[1]{\begin{equation}\label{#1}}
 \newcommand{\ee}{\end{equation}}

 \newcommand{\bl}[1]{\begin{lemma}\label{#1}}
 \newcommand{\el}{\end{lemma}}

 \newcommand{\br}[1]{\begin{remark}\label{#1}}
 \newcommand{\er}{\end{remark}}

 \newcommand{\bt}[1]{\begin{theorem}\label{#1}}
 \newcommand{\et}{\end{theorem}}

 \newcommand{\bd}[1]{\begin{definition}\label{#1}}
 \newcommand{\ed}{\end{definition}}

 \newcommand{\bcl}[1]{\begin{claim}\label{#1}}
 \newcommand{\ecl}{\end{claim}}

 \newcommand{\bfact}[1]{\begin{fact}\label{#1}}
 \newcommand{\efact}{\end{fact}}

 \newcommand{\bp}[1]{\begin{proposition}\label{#1}}
 \newcommand{\ep}{\end{proposition}}

 \newcommand{\bc}[1]{\begin{corollary}\label{#1}}
 \newcommand{\ec}{\end{corollary}}

 \newcommand{\bcj}[1]{\begin{conjecture}\label{#1}}
 \newcommand{\ecj}{\end{conjecture}}

 \newcommand{\bpr}{\begin{proof}}
 \newcommand{\epr}{\end{proof}}

 \newcommand{\bprprop}[1]{\begin{proofof}{\it\ref{#1}}.\,\,}
 \newcommand{\eprprop}{\end{proofof}}

 \newcommand{\bi}{\begin{itemize}}
 \newcommand{\ei}{\end{itemize}}

 \newcommand{\ben}{\begin{enumerate}}
 \newcommand{\een}{\end{enumerate}}

 \newenvironment{proof}{\noindent {\em Proof}.\,\,}{\hspace*{\fill}$\halmos$\medskip}
 \newenvironment{proofof}{\noindent {\em Proof of Proposition\,\,}}{\hspace*{\fill}$\halmos$\medskip}
 \newcommand{\halmos}{\rule{1ex}{1.4ex}}

 \newcommand{\one}{{\mathchoice {1\mskip-4mu\mathrm l}
         {1\mskip-4mu\mathrm l}
         {1\mskip-4.5mu\mathrm l}
         {1\mskip-5mu\mathrm l}}}

\def \Z {{\mathbb Z}}

\def \ra {\rightarrow}
\def \ba {\begin{array}}
\def \ea {\end{array}}

\def \TV {{\hbox{\tiny\rm TV}}}
\def \capac {{\rm capacity}}
\def \Sign {{\rm sgn}}
\def \sgn {{\rm sgn}}

\def \cF {{\mathcal F}}
\def \cG {{\mathcal G}}

\def \cJ {{\mathcal J}}

\def\one{\rlap{\mbox{\small\rm 1}}\kern.15em 1}

\begin{document}
\title{Ergodic behaviour of ``signed voter models''}

\author{\renewcommand{\thefootnote}{\arabic{footnote}}
E.\ Andjel
\footnotemark[1]
\\
\renewcommand{\thefootnote}{\arabic{footnote}}
G.\ Maillard
\footnotemark[2]
\\
\renewcommand{\thefootnote}{\arabic{footnote}}
T.S. Mountford
\footnotemark[2]
}

\date{}

\footnotetext[1]{CMI, Universit\'e de Provence,
13453 Marseille cedex 13, France.
{\sl enrique.andjel@cmi.univ-mrs.fr}
}
\footnotetext[2]{Institut de Math\'ematiques, \'Ecole Polytechnique F\'ed\'erale,
Station 8, 1015 Lausanne, Switzerland.\\
{\sl gregory.maillard@epfl.ch, thomas.mountford@epfl.ch}
}

\maketitle

\begin{abstract}   We consider some questions raised by the recent paper of Gantert, L\"owe and Steif (2005)
concerning ``signed'' voter models on locally finite graphs.  These are voter model like processes with the
difference that the edges are considered to be either positive or negative. If an edge between a site $x$
and a site $y$ is negative (respectively positive) the site $y$ will contribute towards the flip rate of $x$
if and only if the two current spin values are equal (respectively opposed).

\vskip 1truecm
\noindent
{\it MSC} 2000. Primary 60K35, 82C22; Secondary 60G50, 60G60, 60J10.\\
{\it Key words and phrases.} Particle system, voter model, random walk, coupling.\\
{\it Acknowledgment.} The research of EA was partially supported by the European Science Foundation
programme: Phase Transitions and Fluctuation Phenomena. The research of GM and TSM is partially
supported by the SNSF, grants $\#\, 200021-107425$ and $\#\, 200021-107475/1$.
\end{abstract}
\vskip 1 truein
\eject


\begin{section}{Introduction}
\label{S1}

This work arises from questions raised in the recent article by Gantert, L\"owe and Steif, \cite{Gls}.
In this paper we consider voter model like processes called ``signed'' voter models.  For such a process
we suppose given a locally finite graph $G = (V,E)$ and a function
$s\colon E \rightarrow \{-1,1\}$.  Our model $(\eta(t)\colon t \geq 0)$ will simply be a spin system
on $\{-1,1\}^V$ with operator
\be{SVMgendef}
\Omega f( \eta) = \sum_{x \in V} \big(f( \eta^x)- f( \eta)\big) \frac{1}{d(x)}\sum_{y\colon \{x,y\}\in E}
1\{\eta(x)\eta(y) \ne s(\{x,y\})\}.
\ee
Here the usual spins, $0$ and $1$ are replaced by $-1$ and $1$ purely for the resulting notational
simplicity. As usual $d(x)$ is the degree of vertex $x$ and configuration $\eta^x$ is simply the
element of $\{-1,1\}^V$ with spins equal to those of $\eta$ except at site $x$.  From now on we will
abuse notation and write $s(x,y)$ for $s(\{x,y\})$; we will call this the sign of edge $\{x,y\}$.
This can be seen as a generalization of the classical voter model (see e.g. \cite{Liggett85},
\cite{D1}) in that if the function $s$ is identically $1$ (or equivalently if all signs are positive)
then the corresponding process is the voter model.

\bd{drf1}
A nearest neighbour path $(\gamma(s)\colon 0 \leq s \leq t)$ having finitely many jumps at times
$0 \leq t_1 \leq t_2 \leq \cdots t_n \leq t$ is said to be \emph{even} or \emph{positive}  if the
number of $1 \leq i \leq n$ so that $s(\gamma(t_{i^-}), \gamma(t_i)) = -1$  is even.  Otherwise the
path is said to be \emph{odd} or \emph{negative}.  If it is positive we write $\sgn(\gamma) = 1$
otherwise $\sgn(\gamma) = -1$.
\ed

As with the voter model the easiest and most natural way to realize the voter model is via a Harris
construction: we introduce for each ordered pair $(x,y)$ with an edge between them a Poisson process,
$N^{x,y}$, of rate $1/d(x)$ with all Poisson processes being independent.  The process is built
by stipulating that at times $t \in N^{x,y}$, the spin at $x$ becomes equal to $s(x,y)\eta_t(y)$.
A.s.\ no two distinct Poisson processes have common points so the rule is unambiguous.  It can easily
be checked that with probability one this rule specifies $\eta_t( x)$ for all $t$ and $x$ just as in
the classical voter model (see \cite{D1}).  The Markovian nature is simply inherited from that of the
system of Poisson processes.  It is then easily seen that this is indeed the desired process. As with
the voter model, duality plays the dominant role in understanding the ``signed'' voter model. For fixed
$t \geq 0 $ and $x \in V$ we define the random walk on $G$, $X^{x,t}=(X^{x,t}(s)\colon 0 \leq s\leq t)$
by the recipe: $X^{x,t}(0)= x$, the random walk jumps from $y$ to $z$ at time $s\in [0,t]$ if
immediately before time $s$ it was at site $y$ and $t-s \in N^{y,z}$.  As in \cite{D1}, we recover
$\eta_t(x)$ via the identity
\be{dual1}
\eta_t(x) = \eta_0\big(X^{x,t}(t)\big)\, \sgn\big( X^{x,t}\big).
\ee
It should be noted that for fixed $t$ the random walks $X^{x,t}(\cdot) $ and $X^{y,t}(\cdot) $ are
coalescing. If the two paths meet for the first time at $s_o \in [0,t] $, then irrespective of $\eta_0$
we have
\be{dual2}
\eta_t(x) \eta_t(y) = \sgn\big(\gamma^{x,y,s_o}\big),
\ee
where $\gamma^{x,y,s_o}\colon [0,2s_o] \rightarrow V$ is the concatenation of the path
$(X^{x,t}(s)\colon 0 \leq s \leq s_o )$ with the path $(X^{y,t}(s_o-s)\colon 0 \leq s \leq s_o )$.
For more discussion of the dual see the next section.

As written above this article is written to address questions raised by \cite{Gls}; it also follows for
instance the article of \cite{Sa} which addresses signed voter models on the integer lattice where the
signs are assigned to the edges in i.i.d.\ fashion.  See \cite{Gls} for a fuller bibliography.

A major preoccupation of \cite{Gls} was with {\it unsatisfied cycles} that are defined as follows.
\bd{def2}
Unsatisfied cycles are nearest neighbour cycles in $G$ whose sign is negative.
\ed

Such cycles are important since in their absence the vertices can be divided into a ``positive'' set,
$V_+$ and a ``negative'' set, $V_-$ so that the process $(\eta^\prime_t\colon t \geq 0)$ for
$\eta^\prime_t (x)=\eta_t(x) $ for $x \in V_+$, $\eta^\prime_t (x)=- \eta_t(x) $ for $x \in V_-$ is a
classical voter model.  Equally, the presence of unsatisfied cycles precludes the existence of fixed
configurations $\eta$ for which the total flip rate is zero (see \cite{Gls}, Section 2 for details).
For the classical voter model the configurations $\underline 1$ of all $1$s and $- \underline{1}$ of
all $-1$s are fixed in this sense and so the voter model is never ergodic in the sense of
\cite{Liggett85}, i.e., there exists a unique equilibrium $\mu $ and for every initial $\eta_0$,
$\eta_t $ converges in distribution to $\mu $ as $t$ tends to infinity. In the case of ``signed''
voter models ergodicity in this sense is a real possibility. A simple criterion for ergodicity was
the existence of unsatisfied cycles and the recurrence of the associated simple random walk, see
Theorem 1.1 of \cite{Gls}. The question of whether, for these processes, if there existed a unique
equilibrium the process must necessarily be ergodic was raised in \cite{Gls}.  In fact this holds and
can be seen to be a consequence of Matloff's lemma (Lemma 3.1 in \cite{mat77}), see also Lemma V.1.26
of \cite{Liggett85}.
\bt{thm1}
If the ``signed'' voter model  has a unique equilibrium, then the ``signed''
voter model is ergodic.
\et

Another question we are fully able to resolve is the second open question listed in \cite{Gls}:

\bd{def3} For a path $\gamma = (\gamma(s)\colon 0 \leq s \leq t) $ and $0\leq s_1\leq t_1\leq t$,
$\gamma^{s_1,t_1}$ signifies the path $(\gamma(s) \colon s_1 \leq s \leq t_1)$. If $s_1 = 0$ we
write $\gamma^{t_1}$ instead of $\gamma^{s_1,t_1}$. For a path $ \gamma = (\gamma(s)\colon s \geq 0)$
on $V$, we say that $\gamma$ \emph{traverses infinitely many unsatisfied cycles} if there exists
sequences $(s_i)_{i\geq 1}$ and $(t_i)_{i\geq 1}$ tending to infinity so that $\gamma^{s_i,t_i} $ are
unsatisfied cycles.
\ed

\bt{thm2}
For the graph $\Z ^3$ with usual edge set and any sign assignation, $s$, either the process is not ergodic
or a random walk must a.s.\ traverse infinitely many unsatisfied cycles.
\et

By Proposition \ref{prop44} below the two statements in Theorem \ref{thm2} are exclusive.
The peculiarity of this result is highlighted by the next result
\bt{thm3}
For the graph $\Z ^d, \ d \geq 4$ there are sign functions $s$ on the edge set so that the associated
voter model is ergodic but the random walk must a.s.\ traverse only finitely many unsatisfied cycles.
\et

The Theorem 1.1 of \cite{Gls} shows that in dimensions $1$ and $2$, if there is an unsatisfied cycle
then necessarily the associated ``signed'' voter model is ergodic so the above results are in a sense
definitive. We finally consider another raised question (\cite{Gls}, question one). Proposition 1.2
of this work gives a useful robust criterion for there to exist multiple equilibria for a signed voter
model: there exists a subset $W\subset V$ that satisfies firstly that with positive probability a
random walk (starting from an appropriate site) will never leave $W$ and secondly that $W$, with
inherited edge set, has no unsatisfied cycles. The question raised was whether this criterion was in
fact necessary as well as sufficient.
\bp{prop1}
In general if there are multiple equilibria, it does not follow that we can find a region $W \subset V$
on which the inherited graph has no unsatisfied cycle and for which the random walk will with strictly
positive probability never leave.
\ep

But, under a natural condition, the result is in fact true.
\bp{1.2bounded}
If the graph $G=(V,E)$ is of bounded degree and the sign function is such that there are multiple
equilibria, then we can find a region $W \subset V$ on which the inherited graph has no unsatisfied
cycle and for which the random walk will with strictly positive probability never leave.
\ep

Finally in the last section we show
\bp{prop44}
If the random walk on $G=(V, E), (X(t)\colon t \geq 0)$ satisfies with probability 1, $X(\cdot)$ traverses
infinitely many unsatisfied cycles then the signed voter model is ergodic.
\ep

An important tool we will use is the fact that for two Markov chains on a state space $S$ where the jump
rates satisfy
\be{couphyp}
\sup_{x \in S} q(x,x)<\infty ,
\ee
has a ``time shift'' coupling.  By this we mean that
\bl{couplem}
Under condition (\ref{couphyp}), and given $T < \infty $ and $ \epsilon > 0$, there exists
a finite $t_0$ so that for any $s\in [0,T] $ and any $x\in S$, two realizations of the Markov chain starting
at $x$, $(X(t)\colon t\geq 0)$ and $(X^\prime(t)\colon t \geq 0)$ may be coupled so that with probability at
least $1 -\epsilon$
\begin{itemize}
\item[(a)] for all $t \geq t_0$, $X(t)=X^\prime(t+s)$ and
\item[(b)] the sequence of sites visited (allowing repeat visits) by the process $X(\cdot)$ up to time
$t_0$ is equal to that for $X^\prime(\cdot)$ up to time $t_0+s$.
\end{itemize}
\el
(Remark in particular that $\sgn((X^\prime)^{t+s})=\sgn(X^t)$ $\forall\, t\geq t_0$.)

The rest of the paper is organized as follows: Sections \ref{S2}, \ref{S3} and \ref{S4} are respectively
devoted to the proofs of Theorem \ref{thm1}, \ref{thm2} and \ref{thm3}, and Sections \ref{S5} and \ref{S6}
to the proofs of Propositions \ref{prop1}, \ref{1.2bounded} and \ref{prop44}.
\end{section}


\begin{section}{Proof of Theorem \ref{thm1}}
\label{S2}

The following proof for Theorem \ref{thm1} is really just a transcription of Lemma V.1.26 of
\cite{Liggett85}. It is included for completeness. It rests on a property of the dual for the
signed voter model, which we now describe in detail.

We suppose, as usual, a given Harris system for generating signed voter models $(\eta_t\colon t \geq 0)$
from a given initial configuration $\eta_0$.  That is a collection of independent Poisson processes
$N^{x,y}$ of rate $d(x)$ for ordered neighbour pairs $(x,y)$. Given an initial configuration $\eta_0$,
a time $t\geq 0$, an integer $r$ and $r$ points in vertex set $V$, $x_1,x_2,\cdots,x_r$, the values of
$(\eta_t(x_1),\eta_t(x_2),\cdots,\eta_t(x_r))$ are determined by the dual process
\be{S2.3}
\underline{X}^t(u) =
\Big(\big(X^{t,1}(u),i^{t,1}(u)\big), \big(X^{t,2}(u),i^{t,2}(u)\big),\cdots,\big(X^{t,r}(u),i^{t,r}(u)\big)\Big),
\ee
where $X^{t,j}(u) \in V$, $i^{t,j}(u) \in \{-1,1\}$ for all $u\in [0,t]$. The process (piecewise constant)
evolves as follows: $\underline{X}^t(\cdot)$ jumps at time $u\in [0,t]$ if and only if there exists $j\leq r$
so that $t-u \in N^{X^{t,j}(u^-),z}$ for some $z$ neighbouring $X^{t,j}(u^-)$.  This being the case
\begin{itemize}
\item[(i)] for every index $k$ so that $X^{t,k}(u^-) \ne X^{t,j}(u^-)$, there will be no change:
$X^{t,k}(u) = X^{t,k}(u^-)$ and $i^{t,k}(u) = i^{t,k}(u^-)$,
\item[(ii)] for every index $k$ so that $X^{t,k}(u^-) = X^{t,j}(u^-)$, we will have
$X^{t,k}(u) = z$ and $ i^{t,k}(u) = i^{t,k}(u^-)\, s(X^{t,j}(u^-), z)$ with $s$ defined as in
Definition \ref{drf1}.
\end{itemize}
Given this dual one recovers the values $\eta_t(x_k)$ by
\be{S2.5}
\eta_t(x_k) = \eta_0\big(X^{t,k}(t)\big)\, i^{t,k}(t).
\ee
The key point for the proof is that over the interval $[0,t]$ the process $\underline{X}^t$ will evolve as a
Markov chain whose jump rates are bounded and which does not depend on $t$ so that the coupling result mentioned
at the end of the introduction may be applied. That is given integer $r<\infty$ and $\epsilon>0$, uniformly over
all $x_1, x_2,\cdots, x_r$ there exists $t_0$ so that
\be{S2.7}
\Big\|\underline{X}^t(t) - \underline{X}^{t+u}(t+u)\Big\|_\TV
=  \Big\|\underline{X}^{t+u}(t) - \underline{X}^{t+u}(t+u)\Big\|_\TV <  \epsilon
\ee
for all $t\geq t_0 $ and $u\in [0,t]$, where by abuse of notation we identify the random variables with their
law.

We may now turn directly to the proof of Theorem \ref{thm1}. We consider $\eta_0$ as fixed. It is sufficient
to show that all limit points of the distribution of $\eta_{t}$ as $t$ tends to infinity are equilibria. We
suppose that for sequence $\{t_n\}_{n\geq 1}$ tending to infinity
\be{S2.9}
\eta_{t_n} \rightarrow \nu
\quad\text{in law.}
\ee
Let $h$ be a cylinder function depending on, say, the spin values at $x_1, x_2, \cdots, x_r$, i.e.,
$h(\eta) = g(\eta(x_1), \eta (x_2), \cdots, \eta(x_r))$. We have that
\be{S2.11}
< \nu ,h>  = \lim_{n \rightarrow \infty} E^{\eta_0}\big[h(\eta_{t_n})\big]
= \lim_{n \rightarrow \infty} E\big[h^\prime\big(\eta_0, \underline{X}^{t_n}(t_n)\big)\big],
\ee
where by abuse of notation we have
\be{S2.13}
h^\prime\big(\eta_0, \underline{X}^t(t)\big) = g\Big(\eta_0\big(X^{t,1}(t)\big)\, i^{t,1}(t),
\eta_0\big(X^{t,2}(t)\big)\, i^{t,2}(t), \cdots, \eta_0\big(X^{t,r}(t)\big)\, i^{t,r}(t)\Big).
\ee
But equally for any fixed $t$ we have (our signed voter model is easily seen to be a Feller process)
\be{S2.15}
< \nu ,P_t h> = \lim_{n \rightarrow \infty} E^{\eta_0}\big[P_t h(\eta_{t_n})\big],
\ee
where as usual $(P_t)_{t \geq 0}$ denotes the Markov semigroup of our signed voter model. The
quantity inside the limit in the r.h.s.\ of (\ref{S2.15}) can be rewritten as
$E^{\eta_0}[h(\eta_{t_n+t})]$ which in the notation introduced in (\ref{S2.13}) is equal to
\be{S2.17}
E\big[h^\prime\big(\eta_0, \underline{X}^{t_n+t}(t_n+t)\big)\big].
\ee
But, as already noted, as $t_n$ tends to infinity $\|\underline{X}^{t_n}(t_n)-\underline{X}^{t_n+t}(t_n+t)\|_\TV$
tends to zero and so
\be{S2.19}
\lim_{n\to\infty}\Big(E\big[h^\prime\big(\eta_0, \underline{X}^{t_n+t}(t_n+t)\big)\big]
- E\big[h^\prime\big(\eta_0, \underline{X}^{t_n}(t_n)\big)\big]\Big) = 0
\ee
which implies that $< \nu ,h> \,=\, < \nu ,P_t h>$. By the arbitrariness of $t$ and $h$ we must conclude
that measure $\nu $ is an equilibrium but, given our hypotheses that there is a unique equilibrium, we
have established that any limit point $\nu$ must equal this equilibrium. That is we have established ergodicity.

\end{section}


\begin{section}{The integer lattice in three dimensions}
\label{S3}

In this section we consider the signed voter model on $\Z^{3}$ with simple random walk motion.
We address the question of whether the existence of a single equilibrium implies that the simple
random walk must a.s.\ run infinitely many unsatisfied cycles. Given the possibility of adapting
the example of the preceding section to three dimensions we interpret the random walk ``running
infinitely many unsatisfied cycles'' to mean: there exist $s_{i}, t_{i} \uparrow\infty$ with
$s_{i}< t_{i}$ for all $i\geq 1$ so that $B(s_{i})=B(t_{i})$ for all $i\geq 1$ and the path
\be{S3.3}
(B(s)\colon s_{i} \leq s \leq t_{i}):=  B^{s_i,t_i} \text{ is odd}.
\ee
We do not require that the path $B^{s_i,t_i}$ visits each site in the range exactly once, with
the exception of $B(s_i) = B(t_i)$.

Our approach uses the following simple properties of simple random walks found in e.g. Lawler,
\cite{lawler}.
\begin{itemize}
\item[(A)] There exists $k\in (0, \infty)$ so that for a random walk $(X(t) \colon t \geq 0)$ starting
at $X(0)=0$ and any $x\in\partial B(0,n)$
\be{propA}
\frac{1}{kn^{d-1}}
\leq P\Big(X\big(T_{\partial B(0,n)}\big)=x\Big)
\leq \frac{k}{n^{d-1}}
\ee
(see \cite{lawler}, Lemma 1.7.4).
\item[(B)] Harnack principle: for all $\alpha < 1$ there exists $k< \infty$ so that
\be{propB}
\frac{1}{k}
\leq \frac{P^{z}\Big(X\big(T_{\partial B(0,n)}\big)=x\Big)}
{P^{0}\Big(X\big(T_{\partial B(0,n)}\big)=x\Big)}
\leq k
\ee
uniformly over $z\in B(0,\alpha n)$ and $n$ (see \cite{lawler}, Theorem 1.7.6.).
\end{itemize}
Let $C_{r}=\partial B(0,2^r)$, the external boundary, and $B_{r}=B(0,2^r)$. Consider the quantity
\be{S3.4}
H(z)=\sum_{n=1}^{\infty}\sum_{{x \in C_n} \atop {y \in C_{n+1}}}
P\big(X(T_{C_n})=x \,\vert\, X(0)=z\big)\, P\big(X(T_{C_{n+1}})=y \,\vert\, X(0)=x\big)\, N_n^{x,y},
\ee
where for all $x\in C_n$ and $y\in C_{n+1}$
\be{Ndef2}
\begin{aligned}
N_{n}^{x,y}
&= \min\Big\{P^{x}\Big(\text{path }X^{T_{C_{{n+1}}}}
\text{ is even} \,\,\big\vert\,\, X (T_{C_{{n+1}}})= y\Big),\\
&\qquad\qquad P^{x}\Big(\text{path }X^{T_{C_{{n+1}}}}
\text{ is odd} \,\,\big\vert\,\, X(T_{C_{{n+1}}}) = y\Big)\Big\}.
\end{aligned}
\ee
Then, by (\ref{propA}) and (\ref{propB}), the following are clear:
\begin{itemize}
\item[(i)] $H(\cdot)\equiv \infty$ or $H(z)<\infty$ $\forall\,z$;
\item[(ii)] $H(z)<\infty$ if and only if $I<\infty$ with
\be{Idef2}
I = \sum_{n=1}^{\infty}\sum_{{x \in C_n} \atop {y \in C_{n+1}}}
\frac{1}{2^{4n+2}}\, N_{n}^{x,y}.
\ee
\end{itemize}
Furthermore,
\begin{itemize}
\item[(iii)] $I=\infty$ implies that for all random walks a.s.
\be{S3.5}
P\Big(\text{path } X^{T_{C_n}} \text{ is even} \,\,\big\vert\,\, X(T_{C_n})\Big)\to 1/2
\ee
and our voter model is easily seen to be ergodic.
\end{itemize}
Theorem \ref{thm2} will follow from the two following results:
\bp{prop11}
\noindent If $I = \infty$ then a.s.\ the random walk runs infinitely many unsatisfied cycles.
\ep
\bp{prop12}
\noindent If $I < \infty$ then a.s.\ the signed voter model has multiple equilibria.
\ep

\bprprop{prop11}
If $I = \infty$ then one of
\be{S3.7}
\sum_{n = i \textrm{ mod }6} \sum_{{x \in C_n} \atop {y \in C_{n+1}}}
\frac{1}{2^{4n+2}}\, N_{n}^{x,y} = \infty
\ee
for $i = 0,1,2,3,4,5$. Without loss of generality we suppose the first.
The ``mixing'' properties of Brownian motion ensure that then a.s.\
\be{S3.9}
\sum_{n=0 \text{ mod }6} N_{n}^{X(T_{C_n}), X(T_{C_{n+1}})} = \infty
\ee
for any random walk $(X(s)\colon s \geq 0)$. Now we define event $D_{n}$ as
\be{S3.11}
P^{(2^{n},0,0)}\Big(X^\prime \text{ hits } X^{T_{C_{n-2}},T_{C_{n-1}}} \text{ before } C_{n+2}\Big)
\geq C
\ee
for $X^\prime$ an independent random walk, where $C> 0$ is chosen so that for $n$ large
\be{S3.13}
P(D_{n}) > 1/2.
\ee
Define $D^\prime_{n}$ the event
\be{S3.19}
X^{T_{C_{n+1}},T_{C_{n+2}}} \bigcap X^{T_{C_{n-2}},T_{C_{n-1}}}\neq \emptyset.
\ee
By (\ref{propB}) and (\ref{S3.13}), if $D_{n}$ occurs then
\be{S3.15}
P\Big(D^\prime_{n} \,\,\big\vert\,\, \cF_{T_{C_{n-1}}}\Big)> C^\prime,
\ee
for some universal $C^\prime$ not depending on $n$, where $ \{{\mathcal{F}}_t \}_{t \geq 0}$
is the natural filtration for random walk $X(\cdot)$. Now (\ref{propA}) ensures that
\be{S3.17}
\sum_{n= 0 \text{ mod } 6} I_{D^\prime_{n}}\, N^{X(T_{C_{{n}}}), X(T_{C_{{n+1}}})}_ n = \infty
\quad a.s.\
\ee
under conditions given. We now introduce the discrete filtration
\be{S3.20}
\cJ^\prime_{n} = \cF_{T_{C_{6n+2}}}
\quad\text{and}\quad
\cG^\prime_{n} = \sigma \Big(\cJ^\prime_{n}, X^{T_{C_{6n+1}}, T_{C_{6n+2}}}\Big)
\ee
and consider the filtration (over indices $n = 0 \text{ mod 6}$)
\be{S3.21}
\cG^\prime_{1},\cJ^\prime_{1}, \cG^\prime_{2},\cdots,
\cG^\prime_{n},\cJ^\prime_{n}, \cG^\prime_{n+1},\cdots.
\ee
Note that on $D^\prime_{6n} \in \cG^\prime_{n}$ we can define measurably
$t_{n} \in [T_{C_{6n-2}}, T_{C_{6n-1}}]$, $s_{n} \in [T_{C_{6n+1}}, T_{C_{6n+2}}]$
so that $X(t_{n}) = X(s_{n})$. Note that
\be{S3.23}
P\Big(X^{t_n , s_{n}} \text{ is odd} \,\,\big\vert\,\, {\cG}^\prime_{n}\Big)
\geq N_{6n}^{X(T_{C_{6n}}), X(T_{C_{6n+1}})}
\ee
So by (\ref{S3.17}) and L\'evy $0$-$1$ law (see e.g. \cite{dur05}) we have a.s.\ infinitely many
unsatisfied cycles.
\eprprop

\bprprop{prop12}
As before, we denote by $C_r$ the external boundary of $B(0,2^r)$, the Euclidean ball centered
at the origin of radius $2^r$. For $x\in B_{r+1} = B(0,2^{r+1})$, $v\in C_{r+1}$, the law
$P^{x,v,r} $ is the law of the random walk started at $x$ conditioned to exit
$B(0,2^{r+1})$ at $v$. Now for $\alpha<1$ such that $1-\alpha \ll 1$ (and certainly $\leq 1/4$)
and for $x\in C_r$ there are two complementary sets:
\be{S3.25}
S(x,r) = \big\{v\in C_{r+1}\colon N^{x,v}_r < 1-\alpha\big\}
\quad\text{and}\quad
U(x,r) = \big\{v\in C_{r+1}\colon N^{x,v}_r \geq 1-\alpha\big\}.
\ee
For $v \in S(x,r) $ one can speak of a sign of $v$ with respect to $x$: $v$ is even or positive
with respect to $x$ if $P^{x,v,r}(\mbox{path from } x \mbox  { to } v \mbox{ is even}) \geq 3/4$
otherwise $v$ is odd or negative with respect to $x$.  If $v$ is positive with respect to $x$ at
level $r$, we write $\sgn(x,v,r)=1$.  We write $\sgn(x,v,r)=-1$ if $v\in S(x,r)$ (for $\alpha = 3/4$)
but $v$ is not positive with respect to $x$ at level $r$.
For $v\in U(x,r)$ there is (at precision level $1-\alpha$) a reasonable chance of a path from $x$
to $v$ being either even or odd, we write $\sgn(x,v,r)=0$. Therefore, for $u\in C_{r}$ and
$z\in C_{r+1}$,
\be{S3.31}
\Sign(u,z,r) =
\begin{cases}
1  &\text{if } P^{u,z,r}(X\mbox{ is even)} > 3/4;\\
-1 &\text{if } P^{u,z,r}(X\mbox{ is odd} ) > 3/4;\\
0  &\text{otherwise}.
\end{cases}
\ee
We first have
\bl{lem31}
For a random walk $(X(t)\colon t \geq 0 )$ on $\Z^3$ and for any $\alpha < 1$,
under condition $I < \infty $ a.s.\
\be{S3.29}
X(T_{C_{r+1}}) \in S(X(T_{C_{r}}),r)
\ee
for all $r$ sufficiently large.
\el

\bl{lem32}
For any $x \in C_r, w \in C_{r+1}$ with $w \in S(x,r)$, the $P^{x,w,r}$ probability that the path
$X(\cdot)$ satisfies for all $t \leq T_{C_{r+1}} $
\be{S3.33}
\sgn(X^t)\, \Sign(X(t),w,r) = \Sign(x,w,r)
\ee
is at least $4(1-\alpha)$.
\el
\bpr
Suppose without loss of generality that $\Sign(x,w,r)=1$.  Then the $P^{x,w,r}$ probability of event
\be{S3.35}
A = \Big\{\mbox{path } X^{ T_{C_{r+1}}} \mbox{ is odd}\Big\}
\ee
is less than $1- \alpha $.  Consider, with respect to the natural filtration, the c\`adl\`ag martingale
$M_t = E(1_A \,\vert\, \cF_t)$.  By Doob's optional sampling theorem (see e.g. \cite{dur05}) the
probability that this value ever gets above $1/4$ is bounded above by $4(1-\alpha)$. This gives the
result
\epr

The following is a simple consequence of (\ref{propB}).
\bl{lem33}
There exists a universal $c > 0$ so that for any $x \in C_r$ and $w\in C_{r+1} $,
\be{S3.37}
P^{x,w,r}\Big(\capac\big(X^{\tau_r,\sigma _r}\big) > 2^r c\Big) > c\, ,
\ee
where
\be{S3.39}
\tau_r
= \inf\bigg\{t\colon |X(t)| \geq 3\times 2^{r-1} \bigg\}
\quad\text{and}\quad
\sigma_r
= \inf\bigg\{t> \tau_r\colon |X(t)|\geq 7\times 2^{r-2} \mbox{ or }
\leq 5\times 2^{r-2} \bigg\}.
\ee
\el

\bc{cor34}
There exists strictly positive $c$ so that for any $x,y \in C_r $ and $w, v \in C_{r+1} $,
if $X$ is a  $P^{x,w,r}$ motion and $X^ \prime $ is a $P^{y,v,r}$ motion, then with probability
$c$ the conditional probability given $X^\prime$ that $X^{\tau_r,\sigma_r}$ intersects
$(X^\prime)^{\tau_r^\prime, \sigma_r^\prime}$ is at least $c$, where $\tau_r,\sigma_r$ (resp.\
$\tau_r^\prime,\sigma_r^\prime$) are associated to $X$ (resp.\ $X^\prime$).
\ec
\bd{1comp}
We say $\{x,y,v,w\}$ with $x,y\in C_r$ and $v,w\in C_{r+1}$ are $1$-\emph{compatible} if
\be{S3.27}
\sgn(x,v,r)\, \sgn(x,w,r)\, \sgn(y,v,r)\, \sgn(y,w,r)=1.
\ee
\ed
In the following we assume that $\alpha $ has been fixed so large that $240K^2(1- \alpha) < c$ for
$c$ the constant of Corollary \ref{cor34} and $K$ the constant defined in (\ref{S3.47}--\ref{S3.49})
below.

\bl{prop31}
Suppose that $\{x,y,v,w\}$ with $x,y\in C_n $ and $v,w\in C_{n+1}$ are not $1$-compatible and that
for each $u\in\{v,w\}$, $N^{x,u}_{n} < 1 - \alpha$ then for at least one $u\in \{v,w\}$, there exists
some universal constant $K>0$ so that $N^{y,u}_{n} > 3c/(128K)$, where $c$ is the constant
defined in Corollary \ref{cor34}.
\el

\bpr
We suppose without loss of generality that $v$ and $w$ are both positive with respect to $x$ but that
while $w$ is positive with respect to $y$, $v$ is not. By our assumption on the largeness of $\alpha$
we have by Lemma \ref{lem32} and Corollary \ref{cor34}, that there exists a nearest neighbour path
$\gamma(\cdot)$ from $x$ to $w$ on which for all times $s$,
\be{S3.41}
\sgn( \gamma^s)\, \Sign(\gamma(s),w,r) =1.
\ee
Furthermore for $\tau_n^\prime, \sigma_n^\prime$ defined for path $\gamma$, we have
\be{S3.43}
P^{z,u,n}\Big(X^{\tau_n,\sigma_n} \mbox{ hits }
\gamma^{\tau_n^\prime, \sigma_n^\prime} \Big) > c,
\ee
for each $(z,u)\in \{(x,v),(x,w),(y,v),(y,w)\}$. We consider two processes, $(Z^x(t) \colon t\geq 0)$
and $(Z^y(t) \colon t\geq 0)$ starting respectively in $x$ and $y$, running until $C_{n+1} $ is hit
and so that for $u \in \{x,y\}$ the process $(Z^u(t) \colon t\geq 0)$ has law
$1/2\, P^{u,v,n} + 1/2\, P^{u,w,n}$. Then we define the measures $\mu^u(z)$ by
\be{S3.45}
\mu^u(\{z\}) = P\big(Z^u(T_\gamma)=z,\, \tau_n^\prime< T_\gamma < \sigma_n^\prime \big)
\quad\forall\, u\in\{x,y\},\, z\in \gamma^{\tau_n^\prime,\sigma_n^\prime}.
\ee
From facts (\ref{propA}--\ref{propB}), we have that there exists universal $K$ so that
\be{S3.47}
\frac{1}{K}\, \mu^y(\{z\}) \leq
\mu^x(\{z\}) \leq  K\,  \mu^y(\{z\})
\quad\forall\, z \in \gamma^{\tau_n^\prime,\sigma_n^\prime}
\ee
and for either $u$,
\be{S3.49}
P\big(Z^u(T_{C_{n+1}})= v \,\vert\, Z^u(T_\gamma )= z,\,
\tau_n^\prime< T_\gamma < \sigma_n^\prime\big) \in (1/K, 1-1/K)
\quad\forall\, z \in \gamma^{\tau_n^\prime,\sigma_n^\prime}.
\ee
We classify the points in $\gamma $ of size between $5\times 2^{n-2}$ and $7\times 2^{n-2}$
into five sets:
\be{S3.51}
\begin{aligned}
A_{++} &=  \Big\{z\colon P^{z,w,n}\Big(X^{T_{C_{n+1}}} \mbox{ is even}\Big)\geq 3/4,
\,  P^{z,v,n}\Big(X^{T_{C_{n+1}}} \mbox{ is even}\Big) \geq 3/4 \Big\}\\
A_{+-} &=  \Big\{z\colon P^{z,w,n}\Big(X^{T_{C_{n+1}}} \mbox{ is even}\Big)\geq 3/4,
\,  P^{z,v,n}\Big(X^{T_{C_{n+1}}} \mbox{ is odd}\Big) \geq 3/4 \Big\}\\
A_{-+} &=  \Big\{z\colon P^{z,w,n}\Big(X^{T_{C_{n+1}}} \mbox{ is odd}\Big)\geq 3/4,
\,  P^{z,v,n}\Big(X^{T_{C_{n+1}}} \mbox{ is even}\Big) \geq 3/4 \Big\}\\
A_{--} &=  \Big\{z\colon P^{z,w,n}\Big(X^{T_{C_{n+1}}} \mbox{ is odd}\Big)\geq 3/4,
\,  P^{z,v,n}\Big(X^{T_{C_{n+1}}} \mbox{ is odd}\Big) \geq 3/4 \Big\}\\
D      &= \Big\{z\colon P^{z,v,n}\Big(X^{T_{C_{n+1}}} \mbox{ is odd}\Big) \in (1/4, 3/4)\Big\}.
\end{aligned}
\ee
We have by the optimal stopping time reasoning of proof of Lemma \ref{lem32} and our
assumptions on $x$ and $v$ that
\be{S3.53}
\mu^x (D) < 4(1- \alpha ) .
\ee
By (\ref{S3.47}), this implies that
\be{S3.55}
\mu^y (D) < 4(1- \alpha )K .
\ee
We claim that
\be{S3.57}
\mu^x(A_{+-}) < \epsilon = 40K(1- \alpha ).
\ee
To see this suppose the contrary, then we must have either
\be{S3.59}
P\Big(Z^x(T_\gamma) \in A_{+-},\, \tau_n^\prime< T_\gamma< \sigma _n^\prime,\,
(Z^x)^{T_\gamma} \mbox{ is even}\Big)
\geq \epsilon/2
\ee
or
\be{S3.61}
P\Big(Z^x(T_\gamma)\in A_{+-},\, \tau_n^\prime< T_\gamma< \sigma_n^\prime,\,
(Z^x)^{T_\gamma} \mbox{ is odd}\Big)
\geq \epsilon/2.
\ee
In the former case we have via our choice of $K$
\be{S3.63}
P^{x,v,n}\Big(X(T_\gamma)\in A_{+-},\, \tau_n^\prime< T_\gamma< \sigma_n^\prime,\,
X^{T_\gamma} \mbox{ is even}\Big)
\geq \epsilon/(2K)
\ee
and so by the Markov property
\be{S3.65}
P^{x,v,n}\Big(X^{T_{C_{n+1}}} \mbox{ is odd}\Big)
\geq \epsilon /(4K) > 1 - \alpha,
\ee
which contradicts our hypothesis on $x$ and $v$.  Similarly in the other case we are
forced to conclude that $P^{x,v,n}(X^{T_{C_{n+1}}} \mbox{ is odd}) > 1-\alpha$.
Arguing similarly with set $A_{+-}$ replaced by $A_{-+}$, we are able to deduce that
\be{S3.67}
\mu ^x(A_{-+}) < \epsilon.
\ee
The Harnack principle (see (\ref{propB})) now permits us to conclude that
$\mu^y(A_{-+} \cup A_{+-}) < 2K\epsilon $. We thus conclude that either
$\mu^y(A_{++}) \geq (c-2K \epsilon )/2$ or $\mu^y(A_{--}) \geq (c-2K \epsilon )/2$.
Without loss of generality we suppose the former.  Note that our assumptions on the
closeness of $\alpha$ to $1$ ensures that $(c-2K \epsilon )/2 > c/3$.  Then for
identical reasons, either
\be{S3.69}
P \Big(Z^y(T_\gamma)\in A_{++},\, \tau_n^\prime< T_\gamma< \sigma_n^\prime,\,
(Z^x)^{T_\gamma} \mbox{ is even}\Big)
\geq (c-2K \epsilon )/4 > c/6
\ee
or
\be{S3.71}
P\Big(Z^y(T_\gamma)\in A_{++},\, \tau_n^\prime< T_\gamma< \sigma_n^\prime,\,
(Z^x)^{T_\gamma}\mbox{ is even}\Big)
\geq (c-2K \epsilon )/4 > c/6.
\ee
Again without loss of generality we suppose the former.  In this case we have
\be{S3.73}
P^{y,w}\Big(X(T_\gamma)\in A_{++},\, \tau_n^\prime< T_\gamma< \sigma_n^\prime,\,
X^{T_\gamma} \mbox{ is even}\Big)
\geq c/(12K)
\ee
and so
\be{S3.75}
P^{y,w}\Big(X^{T_\gamma} \mbox{ is even}\Big)
\geq c/(32K).
\ee
\epr

Consider two independent random walks $X(t)$ and $Y(t)$ then for any $\alpha<1$ and
any $r$ sufficiently large
\be{S3.77}
N^{X(T_{C_r}), X(T_{C_{r+1}})}_{r}< 1-\alpha
\quad\text{and}\quad
N^{Y(T_{C_r}), Y(T_{C_{r+1}})}_{r}< 1-\alpha,
\ee
that is $\{ X(T^X_{C_{r}}), Y(T^Y_{C_{r}}), X(T^X_{C_{r+1}}), Y(T^Y_{C_{r+1}})\}$ are
$1$-compatible.
\bd{2comp}
We say $\{x,y,z,w\}$ with $x\in C_{r-1}$, $y,z\in C_{r}$ and $w\in C_{r+1}$ are
$2$-\emph{compatible} if
\be{S3.79}
\Sign(x,y,r-1)\, \Sign(x,z,r-1)\, \Sign(y,w,r)\, \Sign(z,w,r) = 1.
\ee
\ed
\bl{prop32}
Under the hypothesis that $I < \infty $, for any two independent random walks
$(X(t) \colon t \geq 0)$ and $(Y(t) \colon t \geq 0)$ with probability one
$\{X(T_{C_{r-1}}),X(T_{C_{r}} ),Y(T_{C_{r}}),X(T_{C_{r+1}})\}$ are $2$-compatible
for all $r$ large.
\el
Here, as before $T_{C_{r}}$ as an argument denotes the stopping time appropriate to the process.

\bpr
We will show that with probability one $\{X(T_{C_{r-1}}),X(T_{C_{r}} ),Y(T_{C_{r}}),X(T_{C_{r+1}})\}$
are $2$-compatible for all $r$ large and even.  The proof for $r$ odd is entirely analogous.
We first observe that under the condition $ I < \infty $, we have a.s.
\be{S3.81}
\sum_{r \text{ even}} N^{X(T_{C_{r-1}}), X(T_{C_{r+1}})}_{r+} < \infty\, ,
\ee
where for all $x\in C_{r-1}$ and $y\in C_{r+1}$
\be{S3.83}
\begin{aligned}
N_{r+}^{x,y}
&= \min\Big\{P^{x}\Big(\text{path }X^{T_{C_{{r+1}}}}
\text{ is even} \,\,\big\vert\,\, X (T_{C_{{r+1}}})= y\Big),\\
&\qquad\qquad P^{x}\Big(\text{path }X^{T_{C_{{r+1}}}}
\text{ is odd} \,\,\big\vert\,\, X(T_{C_{{r+1}}}) = y\Big)\Big\}.
\end{aligned}
\ee
Given (\ref{propA}), we have easily that there exists universal constant $K$ so that the probability
that
\be{S3.85}
N^{X(T_{C_{r-1}}), v}_{r-1} \mbox{ or }  N^{v,X(T_{C_{r+1}})}_{r} > 1/100
\ee
or
\be{S3.87}
\Sign(X(T_{C_{r-1}}), v,r-1)\, \Sign(v,X(T_{C_{r+1}}), r) \ne \Sign(X(T_{C_{r-1}}), X(T_{C_{r+1}}),r+)
\ee
is bounded by $K N^{X(T_{C_{r-1}}), X(T_{C_{r+1}})}_{r+}$, where $\Sign(X(T_{C_{r-1}}), X(T_{C_{r+1}}),r+)$
is given its obvious meaning.  The result now follows from (\ref{propA}) again and L\'evy's $0$-$1$ law.
\epr

Given Lemmas \ref{prop31} and \ref{prop32} we can find a path realization $X=X(t,\omega)$ so that for a.s.\
every random walk path $Y=Y(t,\omega)$ the conclusion of the lemmas hold (here we use the notation
$X(t,\omega)$ to underline the fact that we consider a fixed path of the random walk $(X(s)\colon s\geq 0)$
at time $t$). That is let us pick and fix a ``good'' path $X$ so that for a.s.\ path Y we have that for $r$
large, $\{X(T_{C_{r-1}}), X(T_{C_{r}}), Y(T_{C_{r}}), X(T_{C_{r+1}})\}$ are $2$-compatible and
$\{X(T_{C_{r}}), Y(T_{C_{r}}), Y(T_{C_{r+1}}), X(T_{C_{r+1}})\}$ are $1$-compatible and also such that for
any $\alpha < 1$ eventually $N^{X(T_{C_r}),X(T_{C_{r+1}})}_{r} < 1-\alpha $. We will use this path to
designate sites in $C_r$ as positive of negative: we say that $X_{C_1} $ is a positive site, subsequently if
\be{S3.89}
P^{X(T_{C_{r-1}}), X(T_{C_r}),r}\Big(X^{ T_r} \mbox{ is odd}\Big)
\leq 1/100,
\ee
then $X(T_{C_r})$ has the same sign as $X(T_{C_{r-1}})$.
Given this assignation we now assign signs to arbitrary $y\in C_r$. If
\be{S3.91}
P^{X(T_{C_{r-1}}), y,r}\Big(X^{T_r} \mbox{ is odd}\Big)
\leq 1/100,
\ee
then $y\in C_r$ has the same sign as $X(T_{C_{r-1}})$, otherwise it is the opposite.

\bl{propsign}
With probability one there exists a finite random $r_0$ so that either
\be{S3.93}
\forall r \geq r_0, \qquad \sgn\Big(Y^{T_{C_r}}\Big)\, \sgn\big(Y(T_{C_r})\big) = 1
\ee
or
\be{S3.95}
\forall r \geq r_0, \qquad \sgn\Big(Y^{T_{C_r}}\Big)\, \sgn\big(Y(T_{C_r})\big) = -1.
\ee
\el

\bpr
We first observe that for $r$ large enough all the terms $N^{X(T_{C_{r-1}}),X(T_{C_{r}})}_{r-1}$,
$N^{X(T_{C_{r-1}}),Y(T_{C_{r}})}_{r-1}$ are less than, say, $1/100$.  Furthermore by Lemmas
\ref{prop31} and \ref{prop32} for $r$ large, $1$- and $2$-compatibility give
\be{S3.97}
\begin{aligned}
&\Sign\Big(X(T_{C_{r-1}}),X(T_{C_{r}}) ,r-1\Big)\, \Sign\Big(X(T_{C_{r-1}}),Y(T_{C_{r}}) ,r-1\Big)\\
&\qquad\times\Sign\Big(X(T_{C_{r}}),X(T_{C_{r+1}}),r\Big)\,
\Sign\Big(Y(T_{C_{r}}),X(T_{C_{r+1}}) ,r\Big)=1
\end{aligned}
\ee
and
\be{S3.99}
\begin{aligned}
&\Sign\Big(X(T_{C_{r}}),Y(T_{C_{r+1}}),r\Big)\, \Sign\Big(Y(T_{C_{r}}),Y(T_{C_{r+1}}),r\Big)\\
&\qquad\times\Sign\Big(X(T_{C_{r}}),X(T_{C_{r+1}}),r\Big)\,
\Sign\Big(Y(T_{C_{r}}),X(T_{C_{r+1}}),r\Big)=1.
\end{aligned}
\ee
Therefore their product
\be{S3.101}
\begin{aligned}
&\Sign\Big(X(T_{C_{r-1}}),X(T_{C_{r}}),r-1\Big)\,
\Sign\Big(X(T_{C_{r-1}}),Y(T_{C_{r}}),r-1\Big)\\
&\qquad\times\Sign\Big(X(T_{C_{r}}),Y(T_{C_{r+1}}),r\Big)\,
\Sign\Big(Y(T_{C_{r}}),Y(T_{C_{r+1}}),r\Big)=1.
\end{aligned}
\ee
Using our assumptions, we have
\be{S3.103}
\begin{aligned}
\sgn\big(Y(T_{C_{r+1}})\big) &= \sgn\big(X(T_{C_{r}})\big)\,
\Sign\big(X(T_{C_{r}}),Y(T_{C_{r+1}}),r\big)\\
&= \sgn\big(X(T_{C_{r-1}})\big)\, \Sign\big(X(T_{C_{r-1}}),X(T_{C_{r}}),r-1\big)\,
\Sign\big(X(T_{C_{r}}),Y(T_{C_{r+1}}),r\big)\\
&= \sgn\big(X(T_{C_{r-1}})\big)\, \Sign\big(X(T_{C_{r-1}}),X(T_{C_{r}}),r-1\big)\,
\Sign\big(X(T_{C_{r}}),Y(T_{C_{r+1}}),r\big)\\
&\qquad\times \Sign\big(X(T_{C_{r-1}}),Y(T_{C_{r}}),r-1\big)^2\\
&= \sgn\big(X(T_{C_{r-1}})\big)\, \Sign\big(X(T_{C_{r-1}}),Y(T_{C_{r}}),r-1\big)\,
\Sign\big(X(T_{C_{r-1}}),X(T_{C_{r}}),r-1\big)\\
&\qquad\times
\Sign\big(X(T_{C_{r}}),Y(T_{C_{r+1}}),r\big)\, \Sign\big(X(T_{C_{r-1}}),Y(T_{C_{r}}),r-1\big)\\
&= \sgn\big(Y(T_{C_{r}})\big)\, \Sign\big(X(T_{C_{r-1}}),X(T_{C_{r}}),r-1\big)\,
\Sign\big(X(T_{C_{r}}),Y(T_{C_{r+1}}),r\big)\\
&\qquad\times\Sign\big(X(T_{C_{r-1}}),Y(T_{C_{r}}),r-1\big).
\end{aligned}
\ee
Therefore, combining (\ref{S3.101}) and (\ref{S3.103}), we get for all $r$ large
\be{S3.104}
\sgn(Y(T_{C_{r+1}}))=\sgn(Y(T_{C_{r}}))\,\Sign(Y(T_{C_{r}}),Y(T_{C_{r+1}}),r).
\ee
Now, conditional upon to $Y(T_{C_{r}}), Y(T_{C_{r+1}})$, the probability that
\be{S3.105}
\sgn\big(Y^{T_{C_r}}\big)\, \sgn\big(Y(T_{C_r})\big)
\ne \sgn\big(Y^{T_{C_r}}\big)\, \sgn\big(Y(T_{C_{r+1}})\big)
\ee
is simply $N^{X(T_{C_{r}}),X(T_{C_{r+1}})}_{r}$.  Hence the result follows by Lemma \ref{lem31}.
\epr

Define the function
\be{S3.107}
\begin{aligned}
h(x) &= P^x\Big(\mbox{for all } r \mbox{ large } \sgn\Big(Y^{T_{C_r}}\Big)\,
\sgn\big(Y(T_{C_r})\big) = 1 \Big)\\
&\quad - P^x\Big(\mbox{for all } r \mbox{ large } \sgn\Big(Y^{T_{C_r}}\Big)\,
\sgn\big(Y(T_{C_r})\big) = -1 \Big)
\end{aligned}
\ee
and the product measures $\mu_{\pm}$ by $\mu_+(\{ \eta \colon \eta(x) = 1\}) = (1+h(x))/2$,
$\mu_-(\{ \eta \colon \eta(x) = 1\}) = (1-h(x))/2$. We have by L\'evy's $0$-$1$ law and the Markov
property that with probability $1$ $\lim_{t \rightarrow \infty} |h(Y(t))|$ exists and equals 1. So
there exists $x\in \Z^3$ for which $|h(x)|$ is arbitrarily close to 1 and in particular for which
$h(x)\ne 0$. But in this case we have for all $t$ by duality and the Markov property that
\be{S3.109}
P_t\mu_{\pm}\big(\{\eta\colon \eta(x) = 1\}\big) = (1\pm h(x))/2.
\ee
Then using a similar argument as in \cite{Gls} (Section 7, Proof of Proposition 1.2), this implies
non-uniqueness of equilibria.
\eprprop

\end{section}


\begin{section}{The integer lattice in dimensions four and higher}
\label{S4}

We show Theorem \ref{thm3} in this section.  For notational convenience we give the proof
for four dimensions but the proof is easily seen to hold in all dimensions.

Our purpose is to choose a sequence of integer scales $R_n$ so that $R_{n+1}/R_n$ tends to
infinity sufficiently rapidly. Then we will give sign $+1$ to all edges except those of the
form $(x,x+e_1)$ for $e_1=(1,0,0,0)$ and $x\in \{R_n\}\times[-R_n,R_n]^3$. The basic idea is to
consider a random walk starting at a site in $[-R_n/2,R_n/2]^4$, say, and run until it hits
$\partial [-4R_n,4R_n]^4$. Even given the initial and final points uncertainty as to the sign
of the random walk will be introduced.

In the first part of this section we argue from invariance principle considerations that if
$R_{n+1}/R_n \geq 2(n+1)/K_{n+2}$ for constants $(K_n)_{n\geq 2}$ small then almost surely a random
walk does not run though infinitely many unsatisfied cycles. Then we argue that if we increase the
requirement to
\be{Rndef}
R_{n+1}/R_n\geq 2(n+1)^2/K_{n+2},
\ee
then we will have ergodicty.

We now undertake the first part of the program. Consider a Brownian motion in 4 dimensions,
$(B(t)\colon t \geq 0)$. Let $V_r^i$, $r\geq 0$ and $i=3,4$ be the cube ${[-r, r]}^{i}$ and
given a process $(Y(t) \colon t \geq 0)$, $T_{2n} = \inf\{t\colon Y(t) \mbox{ leaves } V_{2n}^{4}\}$.
It follows from the a.s.\ nonexistence of double points for 4-dim Brownian motion
(see \cite{DEK}) (and the fact that two dimensional subspaces of $\partial V_1^4$ are polar) that,
with probability 1, there does not exist $t_{1}, t_{2} \leq T_{2n}$ so that $t_{1}< t_{2}$ and
\be{S4.3}
\left.
\begin{array}{c}
(B(t_{1}),B(t_2))\\
\text{or}\\
(B(t_{2}),B(t_1))
\end{array}
\right\} \in \big(\{1\} \times V_1^3\big) \times
\big(\partial V_1^4 \setminus (\{1\} \times V_1^3)\big)
\ee
and
\be{S4.5}
B(t_{3}) = B(t_{4}) \text{ for } t_{3} \leq t_{1}  \leq t_{2} \leq t_4.
\ee
Bearing in mind that $B$ does not hit the intersections of the faces of $\partial V_1^{4}$,
there exists $K_{n}>0$ so that with probability greater than $ 1-1/2{n^{2}}$
\be{S4.7}
K_{n} \leq \inf |B(t_{3})-B(t_{4})|
\ee
for $t_{1}, t_{2}, t_{3} , t_{4}$ as above. Now (possibly reducing $K_{n}$) we can also have that
this is so for Brownian motion starting in $V_{K_{n}}^4$ uniformly over the initial point.
Now let us inductively define $R_{n}$ as follows: $R_{1}$ is such that for a $4$-dimensional random
walk (starting at $0$) $X$, the probability that
\begin{itemize}
\item[(i)] there exist $t_{1} < t_{2} \leq T_{2R_{1}}$ (recall that $T_{2R_1}$ is the leaving time of
$V_{2R_{1}}^4$) so that
$$
\left.
\begin{array}{c}
(X(t_{1}),X(t_2))\\
\text{or}\\
(X(t_{2}),X(t_1))
\end{array}
\right\} \in \Big(\{R_{1}\} \times V_{R_{1}}^3\Big) \times
\Big(\partial V_{R_{1}}^4 \setminus \big(\{R_{1}\} \times V_{R_{1}}^3\big)\Big);
$$
\item[(ii)] there exists $t_{3} \leq t_{1}  \leq t_{2} \leq t_4$ so that
$t_{4} \leq T_{2R_{1}}$ and $|X(t_{3})-X(t_{4})| \leq K_{2} R_{1}$
\end{itemize}
is less that $\leq 3/4$. Such an $R_{1}$ exists by the invariance principle, see e.g.\ \cite{dur05}.
Now, given $R_{j-1}$ take $R_{j} \geq 2j\, R_{j-1}/K_{j+1}$, so that for any random walk
$X(\cdot)$ starting in $V_{K_{j+1} R_{j}}^4$, the probability that
\begin{itemize}
\item[(i)] there exists $t_{1}<t_{2}<T_{2(j+1)R_{j}}$ so that
$$
\left.
\begin{array}{c}
(X(t_{2}),X(t_{1}))\\
\text{or}\\
(X(t_{1}),X(t_{2}))\\
\end{array}
\right\} \in \Big(\{R_{j}\} \times V_{R_{j}}^3\Big) \times
\Big(\partial V_{R_{j}}^4 \setminus \big(\{R_{j}\} \times V_{R_{j}}^3\big)\Big);
$$
\item[(ii)] there exists $t_{3}\leq t_{1}\leq t_{2}\leq t_4$ so that
$t_{4} \leq T_{2(j+1)R_{j}}$ and $|X(t_{3})-X(t_{4})| \leq K_{j+1}R_{j}$
\end{itemize}
is bounded by $1/4{(j+1)}^{2}$. Now take the configuration of $\pm 1$ bonds on
${\Z}^{4}$ as follows: all bonds are +1 except bonds
\be{S4.11}
({x}, {x}+e_{1}) \text{ for }x\in \{R_{j}\} \times V_{R_{j}}^3.
\ee
Then by Borel-Cantelli there exists $j_0< \infty$ such that for all $j\geq j_{0}$
\begin{itemize}
\item[(i)] If we consider the random walk between hitting $V_{R_{j}}^3$ until hitting
$V_{j R_{j}}^3$ there is no $t_{1}<t_{2}$ so that
$$
\left.
\begin{array}{c}
(X(t_{2}),X(t_{1}))\\
\text{or}\\
(X(t_{1}),X(t_{2}))\\
\end{array}
\right\} \in \Big(\{R_{j}\} \times V_{R_{j}}^3\Big) \times
\Big(\partial V_{R_{j}}^4 \setminus \big(\{R_{j}\} \times V_{R_{j}}^3\big)\Big);
$$
there exists $t_{3} \leq t_{1}  \leq t_{2} \leq t_4$ such that $X(t_{3}) = X(t_{4})$.
\item[(ii)] The random walk does not return to $V_{R_{j}}^4$ after hitting $V_{j R_{j}}^4$.
\end{itemize}
This easily implies that $X$ does not run thought infinitely many unsatisfied cycles.

It is easily seen that the Harnack principle (property (B) in section \ref{S3}) yields:
\bl{S4.lem1}
Let $\pi_r(w,\cdot)$ be the harmonic measure for a random walk starting at $w$, at the
boundary of the ball $B(0,r)$. Then
\be{S4.13}
\lim_{m\to\infty}\limsup_{r\to\infty}\sup_{{x,y\in B(0,r)}\atop{z\in\partial B(0,mr)}}
\frac{\pi_{mr}(x,z)}{\pi_{mr}(y,z)}=1.
\ee
\el
Let $(R_n)_{n\geq 1}$ satisfying (\ref{Rndef}) and consider $C_n=\partial B(0,n)$ and $S_n=nR_n$.
\bl{S4.lem2}
There exists $k_1\in(0,1/2)$ so that for all $n$ large enough and all $x\in C_{S_n}$, $y\in C_{S_{n+1}}$
\be{S4.15}
P^x\Big(X^{T_{C_{S_{n+1}}}} \text{is odd } \,\big\vert\, X\big(T_{C_{S_{n+1}}}\big)=y\Big)>k_1
\ee
and
\be{S4.17}
P^x\Big(X^{T_{C_{S_{n+1}}}} \text{is even } \,\big\vert\, X\big(T_{C_{S_{n+1}}}\big)=y\Big)>k_1.
\ee
\el
\bpr
By the invariance principle we have that if $n$ is large, uniformly for each $x\in C_{S_n}$ the
probability of leaving the box $V_{R_{n+1}}^{4}$ for the first time through
$\{R_{n+1}\}\times V_{R_{n+1}/2}^{3}$, then passing to $\partial V_{2R_{n+1}}^{4}$ without
leaving $[2R_{n+1}/3,+\infty)\times V_{2R_{n+1}/3}^{3}$ is greater than $k_2\in (0,1)$
for some universal $k_2$. From here, uniformly over the random hitting point of
$\partial V_{2R_{n+1}}^{4}$, the conditional probability of hitting $\partial B(0,(n+1)R_{n+1}/2)$
before hitting $V_{R_{n+1}}^{4}$ will be greater than $k_3\in (0,1)$ provided $n$ is large.
This follows from the invariance principle and the classical hitting estimates of Lawler (see
properties (A) and (B) of Section \ref{S3}). From property (A) of Section \ref{S3} we have the
existence of a constant $k_4$ so that
\be{S4.19}
\frac{1}{k_4\, S_{n+1}^3}\leq P^w\Big(X\big(T_{C_{S_{n+1}}}\big)=z\Big)\leq \frac{k_4}{S_{n+1}^3}.
\ee
So using
\be{S4.21}
\begin{aligned}
&P^w\Big(X\big(T_{C_{S_{n+1}}}\big)=z, T_{V_{R_{n+1}}^4}>T_{C_{S_{n+1}}}\Big)\\
&\qquad\geq P^w\Big(X\big(T_{C_{S_{n+1}}}\big)=z\Big)
-\sup_{z\in V_{R_{n+1}}^4} P^z\Big(X\big(T_{C_{S_{n+1}}}\big)=z\Big)
P^w\Big(T_{C_{S_{n+1}}}>T_{V_{R_{n+1}}^4}\Big),
\end{aligned}
\ee
we obtain
\be{S4.23}
P^w\Big(X\big(T_{C_{S_{n+1}}}\big)=z, T_{V_{R_{n+1}}^4}>T_{C_{S_{n+1}}}\Big)
\geq \frac{1}{2k_4 S_{n+1}^3}
\ee
for $n$ large uniformly over $w\in\partial B(0,(n+1)R_{n+1}/2)$. Hence for all $x\in C_{S_n}$
and $y\in C_{S_{n+1}}$
\be{S4.25}
P^{x}\Big(X^{T_{C_{S_{n+1}}}} \text{ is odd }\,\big\vert\, X\big(T_{C_{S_{n+1}}}\big)=y\Big)
\geq \frac{k_2 k_3}{2k_4 S_{n+1}^{3}}.
\ee
This given,
\be{S4.27}
P^{x}\Big(X\big(T_{C_{S_{n+1}}}\big)=y\Big)\leq \frac{k_4}{S_{n+1}^{3}}
\ee
gives
\be{S4.29}
P^{x}\Big(X^{T_{C_{S_{n+1}}}} \text{ is odd }\,\big\vert\, X\big(T_{C_{S_{n+1}}}\big)=y\Big)
\geq \frac{k_2 k_3}{2k_4^2}.
\ee
We argue similarly for the second part.
\epr

The following is a simple consequence of Lemmas \ref{S4.lem1} and \ref{S4.lem2}.
\bc{S4.cor1}
There exists $N_0$ so that for all $r\in\Z_+$ and all $x\in C_{S_{N_0}}$, $y\in C_{S_{N_0+r}}$
\be{S4.31}
\bigg|2P^x\bigg(X^{T_{C_{S_{N_0+r}}}} \text{ is odd }\,\Big\vert\,\,
X\Big(T_{C_{S_{N_0+r}}}\Big)=y\bigg)-1\bigg|
\leq (1-k_1)^r
\ee
\ec
We are now ready to complete the proof of Theorem \ref{thm3}.
\bp{S4.prop1}
For $R_n$ and sign functions as previously described, the signed voter model is ergodic.
\ep
\bpr
We need only show that as $t$ tends to infinity the difference in absolute variation of the
measures $\mu_{x,t,-}$ and $\mu_{x,t,+}$ tends to zero for each $x\in\Z^4$, where
$\mu_{x,t,\pm}$ is defined on $\Z^4 $ by
\be{S4.33}
\mu_{x,t,+}(y) \ = \ P^{x}\big(X(t)= y, X^ t \mbox{ is even}\big)
\quad\text{and}\quad
\mu_{x,t,-}(y) \ = \ P^{x}\big(X(t)= y, X^ t \mbox{ is odd}\big).
\ee
However we have by our basic coupling that for any $y\in\Z^4$ and any $T\geq 0$
\be{S4.35}
\lim_{t\to\infty}\sup_{s\in [0,T]}\Big(\big\|\mu_{y,t,+}-\mu_{y,t-s,+}\big\|_{\TV}
+\big\|\mu_{y,t,-}-\mu_{y,t-s,-}\big\|_{\TV}\Big)=0
\ee
(see Lemma \ref{lem45} for a statement and a proof in a more general setting).
We consider $x\in\Z^4$, $r\geq 0$ and $T<t$ fixed. Let $\nu_r(ds,y)$ be the joint law
(under $P^x$) of $(T_{C_{S_r}},X(T_{C_{S_r}}))$. Then, by stong Markov property,
\be{S4.37}
\begin{aligned}
\mu_{x,t,+}&=\int_0^T\, \sum_{y\in C_{S_r}}\bigg[
P^x\Big(X^{T_{C_{S_{r}}}}\text{ is odd }\,\big\vert\,\,
T_{C_{S_r}}=s,\, X\big(T_{C_{S_r}}\big)=y\Big)\mu_{y,t-s,-}\\
&\qquad
+P^x\Big(X^{T_{C_{S_{r}}}}\text{ is even }\,\big\vert\,\,
T_{C_{S_r}}=s,\, X\big(T_{C_{S_r}}\big)=y\Big)\mu_{y,t-s,+}
\bigg]\nu_r(ds,y)+\mu_{x,t,+}^r\,,
\end{aligned}
\ee
where
\be{S4.39}
\mu_{x,t,+}^r(z)=P^x\Big(X(t)=y,\, T_{C_{S_r}}>T,\, X^t \text{ is even}\Big),
\ee
and similarly for $\mu_{x,t,-}^r$. But as $t\to\infty$
\be{S4.41}
\sup_{y\in C_{S_{r+1}}}\sup_{0\leq s\leq T}\big\|\mu_{y,t-s,-}-\mu_{y,t,-}\big\|_{\TV}
\longrightarrow 0
\ee
and similarly for $\mu_{x,t-s,+}$. Thus as $t\to\infty$
\be{S4.43}
\begin{aligned}
\mu_{x,t,+}&=\int_0^T\, \sum_{y\in C_{S_r}}\bigg[
P^x\Big(X^{T_{C_{S_{r}}}}\text{ is odd},\, T_{C_{S_r}}\leq T\,\big\vert\,\,
X\big(T_{C_{S_r}}\big)=y\Big)\mu_{y,t,-}\\
&\qquad
+P^x\Big(X^{T_{C_{S_{r}}}}\text{ is even},\, T_{C_{S_r}}\leq T\,\big\vert\,\,
X\big(T_{C_{S_r}}\big)=y\Big)\mu_{y,t,+}
\bigg]\nu_r(ds,y)+\mu_{x,t,+}^r+o(1)
\end{aligned}
\ee
and similarly for $\mu_{x,t,-}$. Thus
\be{S4.45}
\begin{aligned}
\big\|\mu_{x,t,+}-\mu_{x,t,-}\big\|_{\TV}
&=\sum_{y\in C_{S_r}}\bigg[
P^x\Big(X^{T_{C_{S_{r}}}}\text{ is odd},\, T_{C_{S_r}}\leq T\,\big\vert\,\,
X\big(T_{C_{S_r}}\big)=y\Big)\mu_{y,t,-}\\
&\qquad
+P^x\Big(X^{T_{C_{S_{r}}}}\text{ is even},\, T_{C_{S_r}}\leq T\,\big\vert\,\,
X\big(T_{C_{S_r}}\big)=y\Big)\mu_{y,t,+}
\bigg]\nu_r([0,T],y)\\
&+\big\|\mu_{x,t,+}^T-\mu_{x,t,-}^T\big\|_{\TV}+o(1).
\end{aligned}
\ee
Now let $t\to\infty$ and then $T\to\infty$ to get
\be{S4.47}
\begin{aligned}
&\lim_{T\to\infty}\lim_{t\to\infty}\big\|\mu_{x,t,+}-\mu_{x,t,-}\big\|_{\TV}\\
&\qquad=\sum_{y\in C_{S_r}}\bigg|
P^x\Big(X^{T_{C_{S_{r}}}}\text{ is odd }\,\big\vert\,\,
X\big(T_{C_{S_r}}\big)=y\Big)
-P^x\Big(X^{T_{C_{S_{r}}}}\text{ is even }\,\big\vert\,\,
X\big(T_{C_{S_r}}\big)=y\Big)\bigg|.
\end{aligned}
\ee
Now letting $r\to\infty$ and using Corollary \ref{S4.cor1}, we obtain desired the result.
\epr

\end{section}


\begin{section}{Proof of Propositions \ref{prop1} and \ref{1.2bounded}}
\label{S5}

\noindent Gantert et al. \cite{Gls} ask whether the converse of Proposition 1.2 of their article held.
This stated that if the graph $G=(V,E)$ had the property that there existed $W \subset V, x \in V$ so that
\begin{itemize}
\item[(i)] $P^{x} (T_{W} = \infty) > 0$ where $T_{W} = \inf \{t\colon X(t) \in W^c\} $ and
\item[(ii)] $W$ with its inherited edge set contained no unsatisfied cycles,
\end{itemize}
then ``necessarily'' the signed voter model could not be ergodic. The question was raised at the end of
the paper as to whether a converse existed: can it be that whenever a signed voter model is non ergodic
such a $W$ can be found? We first show this is not the case, but then show that with the additional
hypothesis that the graph is of bounded degree, it is indeed true. We first state without proof (it
follows from \cite{Gls}, Proposition 1.2).
\bp{3.1}
If a.s.\ for all random walk $X=(X(t) \colon t\geq 0)$ on the graph $G$, there exists random $T$ so that
on $[T,\infty)$, $X$ does not traverse a negative edge then the signed voter model has multiple equilibria.
\ep

We will build our counterexample out of a rooted tree with only positive edges by adding a number of
negative edges whose density is so small that the property of multiple equilibria is unchanged.
Consider a rooted tree so that each $i^{th }$ generation has $n_{i}$ ``children'' where $n_{i}\ra\infty$
as $i \rightarrow \infty $ and is always even. We now amend $T$ as follows. We pick strictly increasing
$V_{n} \uparrow \infty$ so that $n_{V_{n}} \geq 2^{n}$. At the $V_{n}^{\text{th}}$ generation we pair up
the vertices so that each vertex of the $V_{n}^{\text{th}}$ generation is paired with a member having the
same father. We add the corresponding edges. For the resulting graph all original edges  are fixed positive
and the extra ``within generation'' edges negative.  Though this new graph has cycles, we retain use of the
words descendants inherited from the original rooted tree. By the Borel-Cantelli lemma and Proposition
\ref{3.1}, the signed voter model has multiple equilibria. Let $W$ be a subject of $V$ with the property
that, with initial point suitably chosen, the probability of a random walk on $G$ ever leaving $W$ is
strictly positive.  Then by L\'evy's $0$-$1$ law (see e.g.\ \cite{dur05}), on the event that the random
walk $(X(t)\colon t \geq 0)$ never leaves $W$ we must have with probability tending to $1$ as $t $ tends
to infinity
\be{S5.3}
X\big(T_{r(t)}\big) \mbox{ is a descendant of }X(t)
\ee
and
\be{S5.5}
\mbox{ both } X\big(T_{r(t)}\big) \mbox{ and  its pair belong to } W,
\ee
where $r(t) $ is the next $V_n$ level below the current level of $X(t)$ and $ T_{r(t)}$ is the hitting
time of this generation. But this must mean that with probability tending to one as $t$ tends to infinity,
the cycle of length $3$ involving the point $X(T_{r(t)})$, its pair and their (common) father is
unsatisfied.

This counterexample is somewhat cheap, the ``real'' question is whether the converse to Proposition 1.2
holds for graphs of bounded degree. We now show Proposition \ref{1.2bounded}.

In the following let $M =\sup_{x \in V} d(x)$. We first consider that there exists an equilibrium $\mu$
so that for some $x\in V$, $\mu(\{\eta \colon \eta (x) = 1\}) \neq 1/2$. Let
\be{S5.7}
\alpha = \sup_{x\in V}\big|\mu(\{\eta \colon \eta(x)=1\}) - \mu(\{\eta \colon \eta(x)=-1\})\big|.
\ee
Without loss of generality we have
\be{S5.9}
\alpha = \sup_{x\in V}\big(\mu(\{\eta \colon \eta(x)=1\}) - \mu(\{\eta \colon \eta(x)=-1\})\big).
\ee
Now we have (see e.g.\ \cite{Liggett85} or \cite{mat77}) for any $x\in V$ and $t\geq 0$
\be{S5.11}
h(x) := \mu (\{\eta \colon \eta(x)=1\})- \mu(\{\eta \colon \eta(x)=-1\})
= E^{x}\big[h(X(t))\, \sgn(X^{t})\big].
\ee
Now fix $\epsilon> 0$ with $\epsilon\ll 1$ and let $x\in \{y\colon h(y)>\alpha-\epsilon\}$.
For $0\leq t \leq T$, where $T$ is fixed, let
\be{S5.13}
M_{t} = E\big[\eta_{T}(x) =1 \,\vert\, \cF_t\big] - E\big[\eta_{T} (x) =-1 \,\vert\, \cF_t\big]
= h\Big(X_{t}^{x,T}\Big)\, \sgn\Big(\big(X^{x,T}\big)^{t}\Big).
\ee
(Here $\cF_t=$ Harris system on interval $[T-t, T]$.) Note that $|M_{t}|\leq \alpha$ for all $0\leq t\leq T$.
Let $\sigma = \inf \{t \geq 0 \colon \vert M_{t} \vert \leq \alpha - 10 \epsilon\}$ then by optional sampling
theorem we have
\be{S5.15}
\alpha-\epsilon < h(x)= E[M _{\sigma \wedge T}]\leq (\alpha-10\epsilon) P(\sigma\leq T) + \alpha P(\sigma> T),
\ee
from which we deduce $P(\sigma\leq T) \leq 1/10$. Now this (and the arbitrariness of $T$) implies
that if $W$ is the component of $\{y\colon |h(y)|\geq \alpha - 10 \epsilon\}$ containing $x$, then
$P^{x}(T_{W} = \infty)\geq 9/10$. We now show that, provided $\epsilon$ is sufficiently small $W$ has
no unsatisfied cycles: Suppose not and let $x_{0}, x_{1},\cdots, x_{r}$ be an unsatisfied cycle in $W$.
The point is that for all $i\in\{0,\cdots,r\}$
\be{S5.17}
h(x_{i}) = \sum_{y \sim x_{i}} \frac{h(y)\, s(x_{i},y)}{d(x_{i})},
\ee
thus
\be{S5.19}
h(x_{i}) = \frac{h(x_{i+1})}{M}\,  s(x_{i},x_{i+1}) + R_{i}\bigg(\frac{M -1}{M}\bigg),
\ee
where $|R_{i}|\leq \alpha$. From which we have for $h(x_{i})>0$
\be{S5.21}
(\alpha - 10\epsilon) \leq \frac{h(x_{i+1})\, s(x_{i}, x_{i+1})}{M} + \frac{M -1}{M} \alpha.
\ee
That is $h(x_{i+1})\, s(x_{i},x_{i+1})\geq \alpha - 10M\epsilon > 0$ if $\epsilon$ is sufficiently small.
Similarly if $h(x_{i})< 0$, then $h(x_{i +1})\, s(x_{i},x_{i+1})<-(\alpha - 10M\epsilon)<0$ for $\epsilon$
sufficiently small. This gives a contradiction.

In the following a signed random walk (on graph $G$) shall be a process $((X(t),i(t))\colon t \geq 0)$ so that
$(X(t)\colon t\geq 0)$ is a random walk and the process $(i(t)\colon t\geq 0)$ takes values on $\{-1,1\}$, starts
at value $1$ and only changes when $X(\cdot)$ changes. At a jump time $t$ for $X(\cdot)$, we have
$i(t)=i(t^-)\, s(X(t^-),X(t))$. We now suppose that there exists multiple equilibria but that each equilibria has
\be{S5.23}
h(x) = \mu(\{\eta \colon \eta(x)=1\}) - \mu(\{\eta \colon \eta(x)=-1\}) \equiv 0,
\ee
that is for all $x\in V$, $\mu(\{\eta \colon \eta(x)=1\}) = 1/2$. Let us denote by $\mu_{0}$ the canonical
equilibrium where under $\mu_{0}$ the spins $\eta(x_{1}),\cdots, \eta (x_{n})$ can be obtained by
\begin{itemize}
\item[A)] running coalescing signed random walks $(X^{x_{j}}(\cdot), i^j(\cdot))$, $1\leq j\leq n$ for
``$\infty$'' to obtain coalesced classes
\be{S5.25}
C_{1} = \Big\{x_{i^{1}(1)},\cdots, x_{i^{1}(r_1)}\Big\},
C_{2} = \Big\{x_{i^{2}(1)},\cdots, x_{i^{2}(r_2)}\Big\},
\cdots,
C_{n} = \Big\{x_{i^{n}(1)},\cdots, x_{i^{n}(r_n)}\Big\}.
\ee
\item[B)] assigning signs to $x_1,\cdots, x_{n}$ so that $\eta(x_i)$ must be compatible with $\eta(x_j)$
if $x_{i},x_{j}$ belong to same cluster but are independent and equiprobable if they belong to distinct clusters.
\end{itemize}
Note: automatically we have $\mu = \mu_{0}$ if two independent random walks on $G =(V,E)$ must almost surely meet.

Thus, summarizing the foregoing, it will be enough to show that if every equilibrium $\mu$ satisfies
$\mu(\{\eta \colon \eta(x)=1\}) = 1/2$ for all $x\in V$, then there is a unique equilibrium, the canonical
measure $\mu_{0}$. Thus we consider the evolution of the dual
$((X^{1,t}(\cdot), i^{1,t}(\cdot)),\cdots,(X^{n,t}(\cdot), i^{n, t}(\cdot)))$ for $x_{1},\cdots,x_{n}$ fixed but
$t$ variable (and ultimately tending to $\infty$). For any fixed $n$ and coalescing random walks $X^i(\cdot)$
starting at $x_i\in V$ for $i\in\{1,\cdots, n\}$ let event $A(T)$ be defined by
\be{S5.27}
A(T) = \Big\{\exists\, s >T,\, j,k\in \{1,\cdots,n\} \colon X^j(T) \ne X^k(T) \mbox{ but } X^j(s) = X^k(s)\Big\}.
\ee
We then have for all $\epsilon >0$, there exists $T$ so that $P^{x_{1},\cdots,x_{n}}(A(T)) < \epsilon$.
From this we see that to show our result it is sufficient to show for all $r$ and for all sequence of $r$-tuples
$(y^{n}_{1}, y^{n}_{2},\cdots, y^{n}_{r})$ so that $\lim_{n\to\infty}P^{y^{n}_{1}, \cdots y^{n}_{r}}(A(0))=0$,
 we have for all $(d_1,\cdots, d_r) \in \{-1, +1\}^r$
\be{S5.29}
\lim_{n\to\infty}\mu\Big(\Big\{\eta (y^{n}_{j}) = d_{j},\, j = 1,\cdots,r\Big\}\Big)
= 2^{-r}.
\ee

Remark that even though we have supposed that with positive probability two independent random walks may
avoid each other for ever, nothing prevents the existence of an integer $r$ so that for all distinct
$y_1,y_2,\cdots,y_r$, $P^{y_1,y_2,\cdots,y_r}(A(0))=0$.

To make our claim we will argue by induction. The result for $r=1$ is simply our hypothesis on the
equilibria of our signed voter model. Suppose now that the result holds for $r-1$ and suppose given
a sequence of $r$-tuples $(y^{n}_{1},y^{n}_{2},\cdots,y^{n}_{r})$ and $(d_1,d_2,\cdots,d_r)\in\{-1,1\}^r$.
As a building block we consider the following measure $\gamma_{t,y^n_2,y^n_3,\cdots, y^n_r}$ on
$\{-1, 1\}^{V}$ given by $\gamma_{t,y^n_2,y^n_3,\cdots,y^n_r}(A) = \mu (P_{t} A \vert L_{t})$, where as
usual $(P_t)_{t \geq 0}$ denotes the semigroup for the signed voter model and $L_{t}$ is event that
independent signed random walks random walks (also independent of the voter model)
$((X^j(s),i^j(s)) \colon s \geq 0)$ beginning at $y^n_j$ with $2\leq j\leq r$, satisfy
$\eta(X^j(t))\, i^j(t)=d_j$. We have by induction that as $n$ tends to infinity, the probability of event
$L_t$ tends to $2^{-(r-1)}$. It follows, just as in the proof of Theorem \ref{thm1} (given time stretching
properties of the duals) that any limit point of duals $\gamma_{t,y^n_2,y^n_3,\cdots,y^n_r}$ as $t\to\infty$
is an equilibrium. In particular we have
\be{S5.31}
\lim_{t \rightarrow \infty}
\gamma_{t,y^n_2,y^n_3,\cdots,y^n_r} \big(\big\{\eta\colon \eta(y^n_1) = d_1\big\}\big) = 1/2.
\ee
This implies that for $\epsilon>0$ and $n$ sufficiently large independent signed random walks random walks
(also independent of the $\eta_0$) $((X^j(s),i^j(s)) \colon s\geq 0)$ beginning at $y^n_j$ with
$j\in\{1,\cdots,r\}$, we have
\be{S5.33}
\big|\mu\big(\big\{\eta\colon \eta(X^j_{t})\, i^j(t) = d_j,\, j=1,\cdots,r\big\}\big) - 2^{-r } \big|
< \epsilon
\ee
for $t$ large enough. But, given our assumptions on the sequence $(y^n_1,y^n_2,\cdots,y^n_r)$, implies that
for $n$ large and then $t$ large enough for coalescing signed random walks $((Y^j(s), i^j(s)) \colon u\geq 0)$
starting at $y^n_j$ with $j \in \{1,\cdots,r\}$, we have
\be{S5.35}
\big|\mu\big(\big\{\eta\colon \eta(Y^j_{t})\, i^j(t) = d_j,\, j=1,\cdots,r\big\}\big) - 2^{-r }\big|
< \epsilon.
\ee
But $\mu $ is an equilibrium and this means
\be{S5.37}
\big|\mu\big(\big\{\eta\colon \eta(y^n_j) = d_j,\, j=1,\cdots,r\big\}\big) - 2^{-r }\big|
< \epsilon.
\ee
The result follows from the arbitrariness of $\epsilon$.

\end{section}


\begin{section}{Proof of Proposition \ref{prop44}}
\label{S6}

To show Proposition \ref{prop44} we will need the following result.  Let, for $x\in V$, $t\geq 0$,
the measures $\mu_{x,t,\pm}$ on $V $ be defined by
\be{S6.3}
\mu_{x,t,+}(y) \ = \ P^{x}\big(X(t)= y, X^ t \mbox{ is even}\big),
\ee
\be{S6.5}
\mu_{x,t,-}(y) \ = \ P^{x}\big(X(t)= y, X^ t \mbox{ is odd}\big).
\ee
\bl{lem45}
For fixed $T\in (0,\infty)$, and $\epsilon > 0$, there exists $T_0< \infty$ so that
uniformly over $s \in [0,T]$, $x \in V$ and $t \geq T_0$
\be{S6.7}
\| \mu_{x,t,+} - \mu_{x,t-s,+}\|_{\TV} \ + \ \| \mu_{x,t,-} - \mu_{x,t-s,-} \|_{\TV} \ < \ \epsilon.
\ee
\el
The proof relies on using the coupling of \cite{mountford} for two continuous time random walks on $V$
starting at $x$, $(X(r)\colon r \geq 0)$ and $(X^\prime(r) \colon r \geq 0)$: the associated discrete
time random walk on $V$ starting from $x$ is chosen to be the same for the two continuous time processes.
To complete the realizations of the continuous time processes it is then just a question of adding the
associated i.i.d.\ exponential random variables giving the resting times at each site:
$\{e_i\}_{i \geq 1}$ for process $(X(r)\colon r \geq 0)$ and $\{e^ \prime _i\}_{i \geq 1}$ for process
$(X^\prime(r)\colon r\geq 0)$. We can chose the two realizations so that for all $n $ large
$\sum_{i=1}^n e_i =\sum_{i=1}^n e_i^\prime +s$. The time for this to occur does not depend on the initial
$x$ and is tight over $s $ in compact intervals.

\bprprop{prop44}
We have always the existence of the equilibrium which is the limit of the distribution
$(\eta_{t}\colon t \geq 0)$ for $(\eta_{0} (x))_{x \in V}$ i.i.d. Bernoulli $(1/2)$
with associated distribution $\mu$.
So we must show that for any initial r and $x_{1}, x_{2}, \cdots, x_{r} \in V$ the distribution
\be{S6.9}
(\eta_{t}(x_{1}),\cdots, \eta_{t} (x_{r})) \text{ converges to that of } \mu.
\ee
That is for any $\eta_0$ the joint law of
\be{S6.11}
\eta_{0}(X^{x_{1},t}(t))\, \sgn\big((X^{x_{1},t})^t\big), \eta_{0}(X^{x_{2},t}(t))\,
\sgn\big((X^{x_{2},t})^t\big), \cdots,\eta_{0} (X^{x_{r},t}(t))\, \sgn\big((X^{x_{r},t})^t\big)
\ee
converges to that of $\eta(x_{1}), \eta(x_{2}), \cdots, \eta(x_{r})$ under $\mu$.
Now the $(X^{x_{i},t}(s) \colon 0 \leq s \leq t)$ are coalescing random walks.
But for fixed $x_{1}, x_{2}, \cdots, x_{r}$ the probability of any further coalescence
of the random walks on interval $[T, t]$ converges to zero as $T \rightarrow \infty$,
uniformly in $t > T$. From this we see that to show the desired ergodicity it is enough
to show for $y_1, \cdots, y_{n}$ fixed in $V$ and $(Z^{y_{i}}(s) \colon s \geq 0)$
independent random walks on $G$, $\{\eta_{0}(Z^{y_{i}}(t))\,\sgn((Z^{y_{i}})^t)\colon 1\leq i\leq n\}$
converges in law as $t \rightarrow \infty$ to that of independent Bernoulli $(1/2)$.

We will use induction on integer $n$. We suppose the desired convergence holds for integer
$n-1$ (which is trivial for $n =1$). It is enough to show that as $t \rightarrow \infty$,
the conditional probability that $\eta_{0}(Z^{y_{1}}(t))\, \sgn((Z^{y_{1}})^{t})=1$ given
$\eta_{0}(Z^{y_{i}} (t))\, \sgn((Z^{y_{i}})^t) \,\,\forall\, 2\leq i\leq n)$ converges to $1/2$
in probability. First fix $\alpha>1/2$.  Fix $\epsilon > 0$, a small strictly positive constant
which will be more fully specified later. Fix $T \gg 1$ to be such that
\be{S6.13}
P \big(Z^{y_{1}}(s) \text{ has not traversed an unsatisfied cycle for } 0 \leq s \leq T\big)
< \epsilon/100.
\ee
We suppose that $t \geq T+T_0$ for $T_0$ given by Lemma \ref{lem45} for this $\epsilon$ and $T$.
Consider the martingale
\be{S6.15}
M_{s} = P \big(\eta_{0} (Z^{y_{1} }(t))\, \sgn((Z^{y_{1}})^t) = 1 \,\vert\,
Z^{y_1}(s),\, \sgn\big((Z^{y_1})^s\big),\, \eta_{0} (Z^{y_{i}} (t))\, \sgn((Z^{y_{i}})^t )
\,\,\forall\, 2\leq i\leq n\big).
\ee
On $\{M_{0}> \alpha\}$, we have, conditional on this initial value, by the optional sampling theorem
from \cite{dur05}
\be{S6.17}
P(T_{\alpha} < T ) \leq  \frac{4(1-\alpha )}{3-2 \alpha}
\ee
for
\be{S6.19}
T_{\alpha} = \inf \bigg\{s\colon M_{s} < \frac{1/2 + \alpha}{2}\bigg\}.
\ee
Thus if $ \epsilon $ is sufficiently small then with strictly positive probability
\begin{itemize}
\item[(i)] for all $ 0 \leq s \leq T, \ M_s > \frac{1/2 + \alpha}{2}$ and
\item[(ii)] there exists $0 \leq s_1 \leq s_2 \leq T$ so that $(X(s) \colon s_1 \leq s \leq s_2)$
traverses an unsatisfied cycle.
\end{itemize}
But by Lemma \ref{lem45} and our assumption on $t$ we have that
$$
\big\| \mu_{X({s_1}),t-s_1,+} - \mu_{X(s_1),t-s_2,+}\big\|_{\TV}
+ \big\| \mu_{X(s_1),t-s_1,-} - \mu_{X(s_1),t-s_2,-}\big\|_{\TV}
<  \epsilon.
$$
This and the fact that $M_{s_1} > (1/2 + \alpha)/2 $ implies that $M_{s_2} < 1-(1/2+\alpha)/2 +2\epsilon$.
But if $\epsilon $ is chosen sufficiently small then this will contradict (i) above.  Thus we have that in
fact for $\alpha > 1/2 $ the conditional probability that $\eta_{0} (Z^{y_{1}}(t))\, \sgn((Z^{y_{1}})^t)=1$
given  $\eta_{0} (Z^{y_{i}}(t))\, \sgn((Z^{y_{i}})^ t \,\,\forall\, 2\leq i\leq n)$  is less than $\alpha$
for $t$ large. We similarly have that it must equally be greater than $1-\alpha$ and we are done.
\eprprop

\end{section}


\vspace{0.5cm}


\begin{thebibliography}{9}


\bibitem{D1} {\sc  Durrett, R.} \textit{Lecture notes on particle systems and percolation},
Belmont, CA: Wadsworth, 1988.

\bibitem{dur05}
{\sc  Durrett, R.} \emph{Probability: theory and examples}, Third edition.
Duxbury Press, Belmont, CA, 2005.

\bibitem{DEK}
{\sc  Dvoretzky, A. Erd{\"o}s, P. and  Kakutani, S.} (1950) Double points of paths of Brownian motion in n-space,
{\em Acta Sci. Math.} {\bf 12},  75--81.

\bibitem{Gls} {\sc Gantert, N, L\"owe, M, and Steif, J.} (2005) The voter model with anti-voter bonds
{\em Ann. Inst. H. Poincar\'e Probab. Statist.} {\bf 41},  767--780.

\bibitem{lawler} {\sc Lawler, G.} (1996) {\it Intersections of Random walks.}  Birh\"auser, Boston.

%
%
\bibitem{Liggett85}  {\sc Liggett, T.M.} (1985) {\em Interacting particle
systems.} Grundlehren der Mathematischen Wissenschaften {\bf 276}.
Springer-Verlag, New York.

\bibitem{mat77} {\sc Matloff, N.S.} (1977) Ergodicity conditions for a dissonant voting model,
{\em Ann. Probab.} {\bf 5}, 371--386.

\bibitem{mountford} {\sc Mountford, T.} (1995) A coupling of infinite particle systems.
{\em J.Math Kyoto Univ.} {\bf 35}, 1, 43--52.

\bibitem{Sa} {\sc Saada, E.} (1995) Un mod\`ele du votant en milieu al\'eatoire
{\em Ann. Inst. H. Poincar\'e Probab. Statist.} {\bf 31},  263--271.

\end{thebibliography}
\end{document}